\documentclass[twoside, 12pt]{amsart}
\usepackage{latexsym}
\usepackage{amssymb,amsmath,amsopn}
\usepackage[dvips]{graphicx}   
\usepackage{color,epsfig}      

\newtheorem{theorem}{Theorem}

\newtheorem{prop}{Proposition}

\newtheorem{defi}{Definition}
\newtheorem{prob}{Problem}

\theoremstyle{definition}
\newtheorem{rem}{Remark}

\DeclareMathOperator{\conv}{conv}

\DeclareMathOperator{\vol}{vol}

\DeclareMathOperator{\bd}{bd}

\newcommand{\K}{\mathcal{K}}
\renewcommand{\S}{\mathbb{S}}
\renewcommand{\Re}{\mathbb R}
\newcommand{\B}{\mathbf B}

\newcommand{\M}{\mathcal{M}}
\newcommand{\area}{\lambda}

\begin{document}

\title[On a normed version of a Rogers-Shephard type problem]{On a normed version of a Rogers-Shephard type problem}
\author[Z. L\'angi]{Zsolt L\'angi}
\address{Zsolt L\'angi\\
        Dept. of Geometry, Budapest University of Technology\\
        Egry J\'ozsef u. 1.\\
        1111 Budapest\\
        Hungary}
\email{zlangi@math.bme.hu}

\thanks{Partially supported by the J\'anos Bolyai Research Scholarship of the Hungarian Academy of Sciences.}


\begin{abstract}
A translation body of a convex body is the convex hull of two of its translates intersecting each other.
In the 1950s, Rogers and Shephard found the extremal values, over the family of $n$-dimensional convex bodies, of the maximal volume
of the translation bodies of a given convex body. In our paper, we introduce a normed version of this problem, and for the planar case,
determine the corresponding quantities for the four types of volumes regularly used in the literature:
Busemann, Holmes-Thompson, and Gromov's mass and mass*. We examine the problem also for higher dimensions, and for centrally symmetric convex bodies.
\end{abstract}

\maketitle

\section{Introduction and preliminaries}\label{sec:introduction}

The volume of the convex hull of convex bodies in the Euclidean $n$-space $\Re^n$ has been in the focus of research since the 1950s.
One of the first results in this area is due to Rogers and Shephard \cite{RS58},
who, besides other cases, investigated this volume for two intersecting translates.
They, for an $n$-dimensional convex body $K$, defined the \emph{translation body} of $K$ as the convex hull
of $K \cup (x+K)$ for some $x \in \Re^n$ satisfying $K \cap (x+K) \neq \emptyset$, and determined the extremal values of the quantity
\begin{equation}\label{eq:deftranslate}
c_{tr}(K) = \frac{1}{\lambda_n(K)} \max \{ \lambda_n(\conv (K \cup (x+K)): (x+K) \cap K \neq \emptyset, x \in \Re^n \}
\end{equation}
over the family of $n$-dimensional convex bodies, where $\lambda_n$ and $\conv$ denote $n$-dimensional Lebesgue measure and convex hull, respectively.
Their conjecture about the convex bodies minimizing $c_{tr}(K)$ remained open for almost fifty years.
A proof of this conjecture, using measures in normes spaces, and another one based on more conventional tools, can be found
in \cite{MM06} and \cite{GHL13}, respectively.

The aim of our paper is to introduce a variant of this problem for normed spaces.
In our investigation we denote the family of $n$-dimensional convex bodies by $\K_n$, and the polar of a set $S$ by $S^\circ$.
If $M \in \K_n$ is symmetric to the origin $o$, the normed space with $M$ as its unit ball is denoted by $\M$.
The Euclidean unit ball with $o$ as its centre is denoted by $\B^n$, and we set $\S^{n-1}= \bd \B^n$ and $v_n = \lambda_n(\B^n)$.
For a point $p \in \Re^n$, $|p|$ denotes its Euclidean norm, and, for $p,q \in \Re^n$, by $[p,q]$ we mean the closed segment with endpoints
$p$ and $q$. For simplicity, we call a plane convex body a \emph{disk}, and denote $2$-dimensional Lebesgue measure by $\lambda$.

Let us recall the well-known fact that any finite dimensional real normed space can be equipped with a Haar measure, and that it is
unique up to multiplication of the standard Lebesgue measure by a scalar.
Depending on the choice of this scalar, one may define more than one version of normed volume.
There are four variants that are regularly used in the literature.
The \emph{Busemann} and \emph{Holmes-Thompson volume} of a set $S$ in an $n$-dimensional normed space
with unit ball $M$, is defined as
\begin{equation}\label{eq:BusHTdef}
\vol^{Bus}_M(S) = \frac{v_n}{\lambda_n(M)} \lambda_n(S) \quad \hbox{and} \quad \vol^{HT}_M(S) = \frac{\lambda_n(M^\circ)}{v_n} \lambda_n(S),
\end{equation}
respectively.
Note that the Busemann volume of the unit ball, and the Holmes-Thompson volume of its polar, are equal to that of a Euclidean unit ball.
For \emph{Gromov's mass}, the scalar is chosen in such a way that the volume of a maximal volume cross-polytope, inscribed in the unit ball $M$
is equal to $\frac{2^n}{n!}$ , and for \emph{Gromov's mass*} (or \emph{Benson's definition of volume}), the volume of a smallest volume parallelotope,
circumscribed about $M$, is equal to $2^n$ (cf. \cite{PT10}).
We denote the two latter quantities by $\vol^{m}_M(S)$ and $\vol^{m*}_M(S)$, respectively.

In the light of the previous paragraph, it is clear that for any fixed normed space, the Euclidean result can be immediately applied.

\begin{theorem}\label{thm:fixednorm}
Let $\M$ be a normed space with volume $\vol_M$.
Then, for any convex body $K \in \K_n$, we have
\[
1+\frac{2v_{n-1}}{v_n} \leq \frac{\max \{ \vol_M(\conv (K \cup (x+K))): (x+K) \cap K \neq \emptyset, x \in \Re^n \} }{\vol_M(K)} \leq 1+n.
\]
\end{theorem}

We observe that there is equality on the left if, and only if $K$ is an ellipsoid (cf. \cite{GHL13}), and on the right if,
and only if $K$ is a \emph{pseudo-double-pyramid} (cf. \cite{RS58}).

In the remaining part we use a different approach.
For any $K \in \K_n$, we say that the \emph{relative norm} of $K$
is the norm with the central symmetral $\frac{1}{2}(K-K)$ of $K$ as its unit ball (cf. \cite{LNT13} or \cite{L11}).
Observe that, up to multiplication by a scalar, the relative norm of $K$ is the unique norm in which $K$ is a body of constant width.
We introduce the following quantities.

\begin{defi}\label{defn:main}
Let $K \in K_n$ and $\M$ be the space with its relative norm.
For $\tau \in \{ Bus, HT, m, m* \}$, let
\begin{equation}\label{eq:maindef}
c_{tr}^\tau(K) = \max \{ \vol_M^\tau(\conv (K \cup (x+K))): (x+K) \cap K \neq \emptyset, x \in \Re^n \}.
\end{equation}
\end{defi}

Note that the quantities in Definition~\ref{defn:main} do not change under affine transformations.
Our aim is to characterize the extremal values of these quantities in the planar case.
To formulate our main result we need to define the following plane convex body.

Consider the square $S_0$ with vertices $\left(\pm \frac{1}{\sqrt{2}}, \pm \frac{1}{\sqrt{2}} \right)$ in a Cartesian coordinate system.
Replace the two horizontal edges of $S_0$ by the corresponding arcs of the ellipse with equation
\[
\frac{x^2}{a^2} + \frac{y^2}{b^2} =1,
\]
where $a=1.61803\ldots$, and $b=\frac{a}{\sqrt{2a^2-1}}$.
Note that the vertices of $S_0$ are points of this ellipse.
Replace the vertical edges of $S_0$ by rotated copies of these elliptic arcs by $\frac{\pi}{2}$.
We denote the plane convex body, obtained in this way and bounded by four congruent elliptic arcs, by $M_0$.
We remark that the value of $a$ is obtained as a root of a transcendent equation, and has the property that
the value of $\lambda(M_0^\circ) \left( \lambda(M_0) + 4 \right)$ is maximal for all possible values of $a > 1$.

Our main result is the following.

\begin{theorem}\label{thm:main}
Let $K \in \K_2$. Then
\begin{itemize}
\item[\ref{thm:main}.1.] we have $2\pi \leq c_{tr}^{Bus}(K) \leq 3\pi$, with equality on the left if, and only if
$K$ is a triangle, and on the right if, and only if $K$ is a parallelogram.
\item[\ref{thm:main}.2.] We have $\frac{18}{\pi} \leq c_{tr}^{HT} (K) \leq 7.81111\ldots$, with equality on the left if , and only if $K$ is a triangle,
and on the right if $K$ is an affine image of $M_0$.
\item[\ref{thm:main}.3.] We have $6 \leq c_{tr}^{m}(K) \leq \pi + 4$, with equality on the left if, and only if $K$ is a (possibly degenerate) convex quadrilateral, and on the right if, and only if $K$ is an ellipse.
\item[\ref{thm:main}.4.] We have $6 \leq c_{tr}^{m^*}(K) \leq 12$, with equality on the left if, and only if $K$ is a triangle, and on the right if, and only if $K$ is a parallelogram.
\end{itemize}
\end{theorem}

It is a natural question to ask for the extremal values of these four quantities over the family of centrally symmetric plane convex bodies.
This question is answered in the next theorem.

\begin{theorem}\label{thm:symmetric}
Let $M \in \K_2$ be $o$-symmetric. Then
\begin{itemize}
\item[\ref{thm:main}.1.] we have $\pi+4 \leq c_{tr}^{Bus}(M) \leq 3$, with equality on the left if, and only if
$M$ is an ellipse, and on the right if, and only if $M$ is a parallelogram.
\item[\ref{thm:main}.2.] We have $\frac{21}{\pi} \leq c_{tr}^{HT} (M) \leq 7.81111\ldots$, with equality on the left if, and only if $M$ is an affine-regular hexagon, and on the right if $M$ is an affine image of $M_0$.
\item[\ref{thm:main}.3.] We have $6 \leq c_{tr}^{m}(M) \leq \pi + 4$, with equality on the left if, and only if $M$ is a parallelogram,
and on the right if, and only if $M$ is an ellipse.
\item[\ref{thm:main}.4.] We have $7 \leq c_{tr}^{m^*}(M) \leq 12$, with equality on the left if, and only if $M$ is an affine-regular hexagon,
and on the right if, and only if $M$ is a parallelogram.
\end{itemize}
\end{theorem}

The proof of Theorem~\ref{thm:symmetric} is a straightforward modification of the proof of Theorem~\ref{thm:main},
and thus, we omit it.

In Section~\ref{sec:2_left}, we prove the left-hand side inequality about Holmes-Thompson area. In Section~\ref{sec:2_right} we deal with the right-hand side inequality regarding it.
In Section~\ref{sec:134} we examine Busemann area, Gromov's mass and its dual. Finally, in Section~\ref{sec:remarks}, we collect our remarks, and propose some open questions.

\section{The proof of the left-hand side inequality in \ref{thm:main}.2} \label{sec:2_left}

Let $K \in \K_2$ and $M = \frac{1}{2}(K-K)$. From (\ref{eq:BusHTdef}) and (\ref{eq:maindef}), one can deduce that
\begin{equation}\label{eq:HTfirst}
c_{tr}^{HT}(K) = \frac{\area(M^\circ)}{\pi} \left( \area(K) + \max \{ d_K(u) w_K(u^\perp) : u \in \S^1\} \right),
\end{equation}
where $d_K(u)$ is the length of a longest chord of $K$ in the direction of $u$, and $w_K(u^\perp)$ is the width of $K$ in the direction perpendicular to $u$ (cf. also the proof of Theorem 1 in \cite{GHL13}).

Observe that for any direction $u$, we have $d_K(u)=d_M(u)$ and $w_K(u)=w_M(u)$, which yields that minimizing $c_{tr}^{Bus}(K)$, over the class of convex disks with a given central symmetral, is equivalent to minimizing $\area(K)$ within this class.
For the special case that $M$ is a Euclidean unit ball, this problem is solved by a theorem of Blaschke \cite{B15} and Lebesgue \cite{L14},
which states that the smallest area convex disks of constant width two are the Reuleaux triangles of width two.
This result was generalized by Chakerian \cite{C66} for normed planes in the following way.

\begin{figure}[here]
\includegraphics[width=0.8\textwidth]{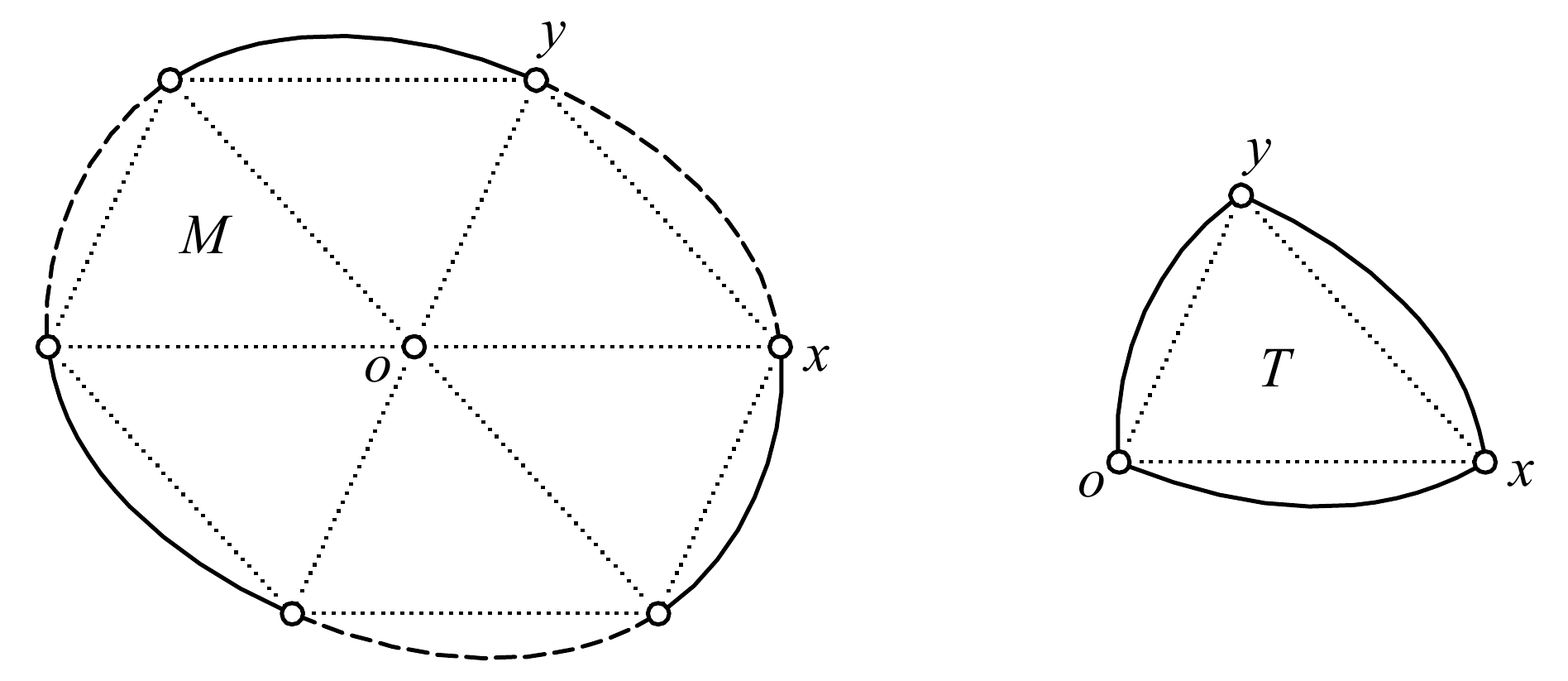}
\caption[]{The construction of Reuleaux triangles in a normed plane}
\label{fig:Reuleaux}
\end{figure}

Let $M \subset \K_2$ be an $o$-symmetric convex disk. Then, for every $x \in \bd M$, there is an affine-regular hexagon, inscribed in $M$, with $x$ as a vertex. Let $y$ be a consecutive vertex of this hexagon. By joining the points $o$, $x$ and $y$ with the corresponding arcs in $\bd M$ we obtain a `triangle' $T$ with three arcs from $\bd M$ as its `sides' (cf. Figure~\ref{fig:Reuleaux}).
These `triangles', and their homothetic copies, are called the \emph{Reuleaux triangles in the norm of $M$}.
Chakerian proved that, given a normed plane $\M$, the area of any convex disk $K$ of constant width two in the norm of $\M$ is minimal for some Reuleaux triangle in the norm. It is not too difficult to see, and was also proven by Chakerian, that the area of such a triangle is equal to $\area(K) = 2 \area(M) - \frac{4}{3}\area(H)$, where $H$ is a largest area affine-regular hexagon inscribed in the unit disk $M$.

Now, assume that $K \in \K_2$ is a minimizer of $c_{tr}^{HT}(K)$ over $\K_2$; by compactness arguments, such a minimizer exists. Then, from Chakerian's result, we obtain that $K$ is a Reuleaux triangle in its relative norm, and that its area is $\area(K) = 2 \area(M) - \frac{4}{3}\area(H)$, where $H$ is a largest area affine-regular hexagon inscribed in $M$.
Now let $P$ be a largest area parallelogram inscribed in $M$.
Then, by (\ref{eq:HTfirst}) and the equality
\[
\max \{ d_K(u) w_K(u^\perp) : u \in \S^1\} = 2 \area(P),
\]
we have
\begin{equation}\label{eq:HTsecond}
c_{tr}^{HT}(K) = \frac{\area(M^\circ)}{\pi} \left( 2 \area(M) - \frac{4}{3}\area(H) + 2 \area(P) \right) .
\end{equation}

It is easy to see that if $K$ is a triangle, then $M$ is an affine-regular hexagon, and vice versa, if $M$ is an affine-regular hexagon,
then the smallest area Reuleaux triangles in its norm are (Euclidean) triangles.
Thus, we only need to show that the quantity in (\ref{eq:HTsecond}) is minimal
if, and only if $M=H$.
Observe that $\area(H) \leq \area(M)$, and hence, it suffices to prove that
\begin{equation}\label{eq:HTthird}
f(M) = \frac{\area(M^\circ) \left( \frac{2}{3} \area(M) + 2 \area(P) \right) }{\pi}
\end{equation}
is minimal if, and only if $M$ is an affine-regular hexagon.

Now we show that if $f(M)$ is minimal for $M$, then its norm is a \emph{Radon norm} (cf. \cite{MSW01} or \cite{PT05}).
Recall that a norm is Radon if, for some affine image $C$ of its unit disk, the polar $C^\circ$ is a rotated copy of $C$ by $\frac{\pi}{2}$;
in this case the boundary of the unit disk is called a \emph{Radon curve}.

Since $f(M)$ is an affine invariant quantity, we may assume that $P$ is a square, with vertices $(\pm 1, 0)$ and $(0,\pm 1)$ in a Cartesian coordinate system.
Note that as $P$ is a largest area inscribed parallelogram, the lines $x=\pm 1$ and $y=\pm 1$ support $M$.
Thus, the arc of $\bd M$ in the first quadrant determines the corresponding part of $\bd M^\circ$.
On the other hand, the maximality of the area of $P$ yields that for any point $p \in \bd M$, the two lines, parallel to the segment $[0,p]$ and at the distance $\frac{1}{|p|}$ from the origin, are either disjoint from $M$ or support it.
Thus, the rotated copy of $M^\circ$ by $\frac{\pi}{2}$ contains $M$, and the two bodies coincide if, and only if $\bd M$ is a Radon curve.

\begin{figure}[here]
\includegraphics[width=0.5\textwidth]{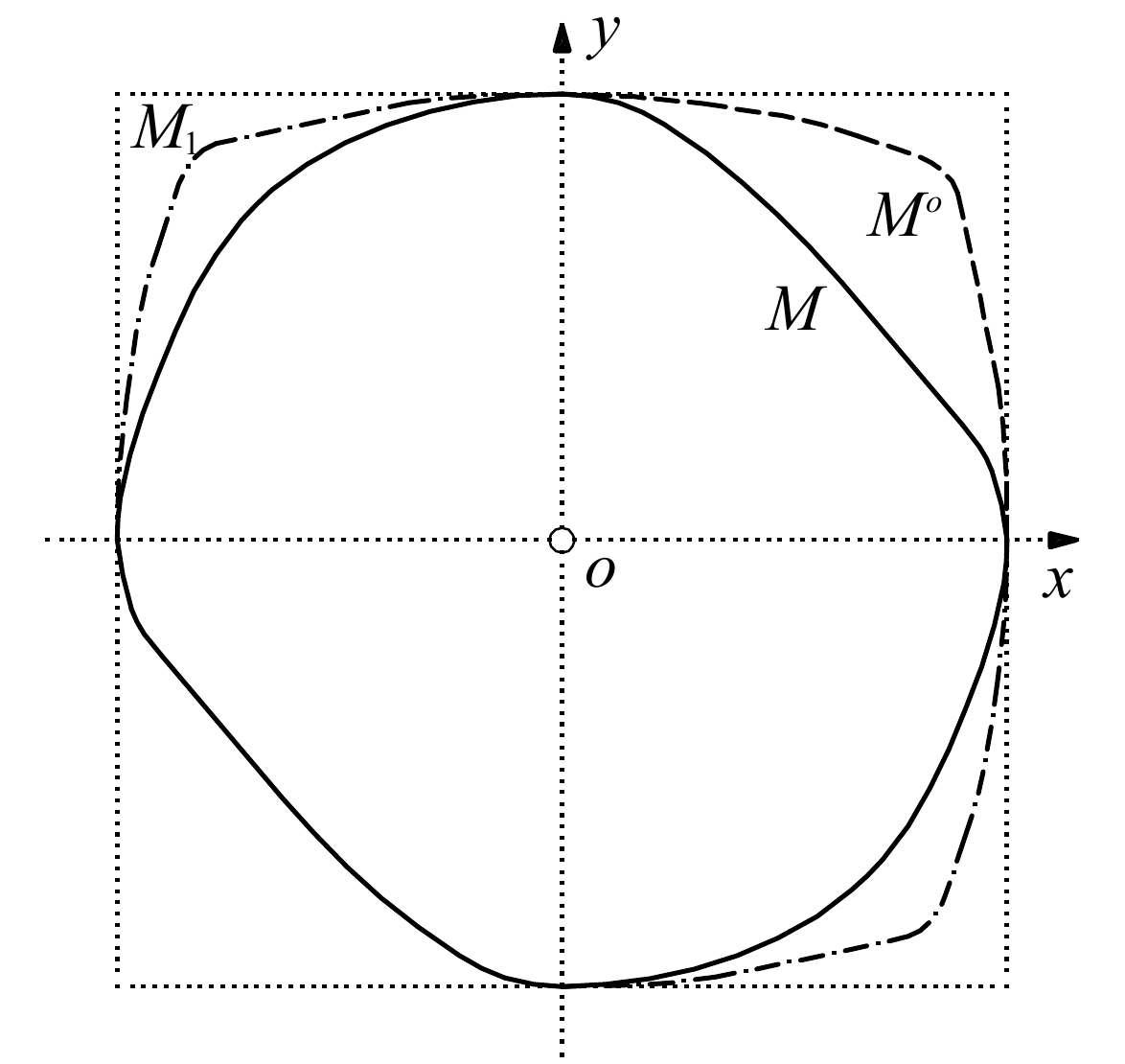}
\caption[]{The extension of $M$ to the unit disk of a Radon norm}
\label{fig:Radon}
\end{figure}

Let $Q_1$ and $Q_2$ denote the parts of $M$ in the first and the second quadrant, respectively.
We define $Q_1^\circ$ and $Q_2^\circ$ similarly for $M^\circ$.
Then $\area(Q_2^\circ) = \area(Q_1) + x_1$ and $\area(Q_1^\circ) = \area(Q_2)+x_2$ for some $0 \leq x_1,x_2 \leq \frac{1}{2}$.
Using this notation, we have $f(M) = \frac{1}{\pi} \left(\area(M)+2x_1+2x_2 \right) \left( \frac{2}{3}\area(M) + 4 \right)$.
Let $M_1$ denote the convex disk obtained by replacing the part of $\bd M$ in the second and fourth quadrants by the rotated copy of the arc of $\bd M^\circ$
in the first quadrant (cf. Figure~\ref{fig:Radon}).
Similarly, let $M_2$ be the disk obtained by replacing the part of $\bd M$ in the other two quadrants by the rotated copy of the arc of $\bd M^\circ$ in the second quadrant.

By our previous observations, we have that $M_1$ and $M_2$ are unit disks of Radon norms, and $M \subset M_1$ and $M \subset M_2$.
On the other hand, the area of a largest area parallelogram inscribed in $M_1$ or $M_2$ is equal to $\lambda(P) = 2$.
Now an elementary computation shows that
\[
f(M_i) = \frac{1}{\pi} \left(\area(M)+2x_{i+1} \right) \left( \frac{2}{3}(\area(M)+2x_{i+1}) + 4 \right) \quad \hbox{for} \, i =1,2,
\]
which, since $0 \leq x_1, x_2 \leq \frac{1}{2}$, yields that
\[
2f(M)-f(M_1)-f(M_2) = \frac{1}{\pi} \left( 8x_1+8x_2-\frac{8}{3}x_1^2-\frac{8}{3}x_2^2 \right) \geq 0,
\]
with equality if, and only if $x_1 = x_2 = 0$. From this, it follows that $f(M) \geq \min \{ f(M_1), f(M_2) \}$, with equality if, and only if
$x_1 = x_2 = 0$ and $M_1 = M_2 = M$. This readily implies that if $f(M)$ is minimal for $M$, then $M$ is the unit disk of a Radon norm.

In the following, we assume that the norm of $M$ is Radon.
Observe that, under our assumption about $P$, we have $\area(M) = \area(M^\circ)$, since $M^\circ$ is a rotated copy of $M$.
On the other hand, since the \emph{volume product} $\area(M) \area(M^\circ)$  of $M$ (cf. e.g. \cite{BM87})
does not change under affine transformations, the definition of Radon norm implies that, in general,
\[
\area(M^\circ) = \frac{4\area(M)}{(\area(P))^2}.
\]
Since $\vol_M^m(M) = \frac{2}{\area(P)} \area(M)$ (cf. the definition in Section~\ref{sec:introduction}, or \cite{PT10}), this yields that
\[
f(M) = \frac{4\area(M)}{(\pi \area(P))^2} \left( \frac{2}{3} \area(M) + 2 \area(P) \right) = \frac{2}{3\pi} \left( \vol_M^m(M) \right)^2 + \frac{2}{\pi} \vol_M^m(M).
\]

Hence, we need to find the minimum of $\vol_M^m(M)$ under the condition that $M$ defines a Radon norm.
This problem was examined in \cite{PT05}, where the authors proved that
for any Radon norm with unit disk $M$, $\vol_M^m(M)$  is at least $3$, with equality if, and only if $M$ is an affine-regular hexagon.
Thus, the left-hand side of \ref{thm:main}.2 immediately follows.

\section{The proof of the right-hand side inequality in \ref{thm:main}.2}\label{sec:2_right}

Assume that $c_{tr}^{HT}(K)$ is maximal for some $K \in \K_2$ and let $M = \frac{1}{2}(K-K)$.
Note that by the Brunn-Minkowski Inequality, we have $\area(K) \leq \area (M)$, with equality if, and only if $K$ is
centrally symmetric. Thus, (\ref{eq:HTfirst}) implies that $K$ is centrally symmetric and, without loss of generality, we may assume that
$K = M$.

Let $P$ be a largest area parallelogram inscribed in $M$.
Since $c_{tr}^{HT}(M)$ is affine invariant, we may assume that $P$ is the square with vertices $(\pm 1, 0)$ and $(0,\pm 1)$ in a Cartesian coordinate system.
Then the lines $x=\pm 1$ and $y=\pm 1$ support $M$.
Let $\sigma$ be a Steiner symmetrization with a symmetry axis of $P$ as its axis, and let $M^*=\sigma(M)$.
Then, clearly, $\area(M^*) = \area(M)$. Observe that $P$ is inscribed in $M^*$ as well, which yields that if $P^*$ is a
maximal area parallelogram inscribed in $M^*$, then $\area(P^*) \geq \area(P)$.
For the Euclidean version of the problem, we have
\begin{equation}\label{eq:ezkell}
c_{tr}(M) = 1 + \frac{2\area(P)}{\area(M)}.
\end{equation}
Then, Theorem 1 of \cite{RS58} yields that $c_{tr}(M)$ does not increase under Steiner symmetrization, which implies
that $\area(P^*) \leq \area(P)$. Thus, we have $\area(P^*) = \area(P)$.

Now we apply a result of Meyer and Pajor \cite{MP90} about the Blascke-Santal\'o Inequality,
who proved that volume product does not decrease under Steiner symmetrizations, which yields that $\area((M^*)^\circ) \geq \area(M^\circ)$.
Thus, since $M$ maximizes $c_{tr}^{HT}(M)$, (\ref{eq:HTfirst}) implies that $\area((M^*)^\circ) = \area(M^\circ)$.
Unfortunately, no geometric condition is known that characterizes the equality case for Steiner symmetrization.
Nevertheless, we may apply another method, used by Saint-Raymond in \cite{SR81}, which he used to characterize the equality case of the Blaschke-Santal\'o Inequality. This method, described also in \cite{W07}, is as follows.

Let $C$ be an $o$-symmetric convex body in $\Re^n$, and let $H$ be the hyperplane with the equation $x_n = 0$.
For any $t \in \Re$, let $C_t$ be the section of $C$ with the hyperplane $\{ x_n = t\}$.
Define $\bar{C}$ as the union of the $(n-1)$-dimensional convex bodies $te_n + \frac{1}{2}(C_t-C_t)$, where $e_n$ is the $n$th coordinate unit vector.
Then we have the following (cf. Lemma 5.3.1 and the proof of Theorem 5.3.2 of \cite{W07}).
\begin{itemize}
\item $\bar{C}$ is an $o$-symmetric convex body.
\item $\lambda_n(\bar{C}) \geq \lambda_n(C)$, with equality if, and only if every $t$-section $C_t$ has a centre of symmetry.
\item $\lambda_n(\bar{C}^\circ) \geq \lambda_n(C^\circ)$.
\item If $\lambda_n(\bar{C}^\circ) \lambda_n(\bar{C})= \lambda_n(C^\circ) \lambda_n(C)$, then the centres of symmetry of the sets $C_t$ lie on a straight line segment.
\end{itemize}

We note that this symmetrization procedure in the plane coincides with the Steiner symmetrization with respect to the second coordinate axis.

Let $L$ be the axis of $\sigma$. Then, since in our case $\area(M^*)=\area(M)$ and $\area((M^*)^\circ)=\area(M^\circ)$,
it follows from the theorem of Saint-Raymond that the midpoints of the chords of $M$, perpendicular to $L$, lie on a straight line segment.
On the other hand, as $\sigma(P)=P$, we have that this segment is contained in $L$.
Thus, $M$ is symmetric to $L$.
Since $L$ was an arbitrary symmetry axis of $P$, we obtain that the symmetry group of $M$ contains that of $P$, and, in particular,
$M$ has a $4$-fold rotational symmetry.

Observe that in this case $M \subseteq \B^2$. Indeed, if for some $p \in M$ we have $|p| > 1$, then, by the $4$-fold rotational symmetry of $M$,
it follows that $M$ contains a square of area greater than $\area(P) =2$, which contradicts our assumption that $P$ is a largest area parallelogram
inscribed in $M$. Since it is easy to check that $c_{tr}^{HT}(M)$ is not maximal if $M = \B^2$, this implies, in particular, that $\lambda(M) < \pi$.
Note that in our case the area of the part of $M$ in each quadrant is equal.

In the next step, we use the following Proposition from \cite{BM13}.

\begin{prop}[B\"or\"oczky, Jr., Makai, Jr.]
Let $Q = \conv \{ o,a,c,b \}$ be a convex deltoid symmetric about the line containing the diagonal $[o,c]$.
Assume that $a,b \in \S^1$ and that the lines containing $[a,c]$ and $[b,c]$ support $\B^2$.
Let $C$ be any $o$-symmetric plane convex body such that $a,b \in \bd C$ and the lines containing $[a,c]$ and $[b,c]$ support $C$,
and set $K=C \cap Q$ and $K^\circ = C^\circ \cap Q$. Let $\lambda (K) = \alpha \leq \lambda (Q \cap \B^2)$ be fixed.
Then $\lambda(K^\circ)$ is maximal, e.g., if $C$ is an $o$-symmetric ellipse $E$ satisfying $\lambda (E \cap Q) = \alpha$.
\end{prop}

Applying this theorem for the part of $M$, say, in the first quadrant, we have that, under our assumption about $P$, $M$ is a convex body bounded by four congruent elliptic arcs, having centres at $o$.
Then it is a matter of computation to verify that $f(M)$ is maximal for a rotated copy of the body $M_0$ described in the introduction.

\section{The proofs of \ref{thm:main}.1, \ref{thm:main}.3 and \ref{thm:main}.4}\label{sec:134}

First, we prove \ref{thm:main}.1.
Observe that for any $K \in \K_2$,
\[
c_{tr}^{Bus}(K) = \frac{\pi}{\area(M)} \left( \area(K) + 2\area(P) \right),
\]
where $M = \frac{1}{2}(K-K)$, and $P$ is a largest area parallelogram inscribed in $M$.
By the result of Chakerian \cite{C66} described in Section~\ref{sec:2_left}, we have that if $K$ minimizes $c_{tr}^{Bus}(K)$ over $\K_2$,
then $K$ is a minimal area Reuleaux triangle in the norm of $M$, and its area is
\begin{equation}\label{eq:Bus_first}
\area(K) = 2 \area(M) - \frac{4}{3}\area(H),
\end{equation}
where $H$ is a largest area affine-regular hexagon inscribed in $M$.
Thus, we may assume, without loss of generality, that
\begin{equation}\label{eq:Bus_second}
c_{tr}^{Bus}(K) = \frac{\pi}{\area(M)} \left( 2\area(M) + 2\area(P) - \frac{4}{3}\area(H) \right).
\end{equation}
Note that in this case $K$ is a (Euclidean) triangle if, and only if $M=H$.

From (\ref{eq:Bus_second}), it readily follows that
\[
c_{tr}^{Bus}(K) = 2\pi + 2\pi \frac{ 3\area(P) - 2 \area(H)}{3 \area(M)}.
\]
Observe that $H$ contains a parallelogram of area $\area(\bar{P})=\frac{2}{3}\area(H)$.
Since $H \subseteq M$, this yields that $\area(P) \geq \frac{2}{3}\area(H)$, with equality if, and only if $M=H$.
This means that $c_{tr}^{Bus}(K) \geq 2\pi$, with equality if, and only if $M=H$, which proves the left-hand side inequality about Busemann area.

Now we prove the right-hand side inequality.
The formula in (\ref{eq:Bus_first}) and the Brunn-Minkowski Inequality shows, like in Section~\ref{sec:2_right}, that
if $c_{tr}^{Bus}(K)$ is maximal over $\K_2$, then $K$ is centrally symmetric.
Thus we may apply Theorem 3 of \cite{RS58} about the maximum of $c_{tr}(K)$, which yields the assertion.

Next, we prove \ref{thm:main}.3.
Let $P$ be a largest area parallelogram inscribed in $M=\frac{1}{2}(K-K)$.
Then we have
\begin{equation}
c_{tr}^{m}(K) = \frac{2\left( \area(K) + 2\area(P) \right)}{\area(P)} = 4 + \frac{2 \area(K)}{\area(P)}.
\end{equation}
Observe that for any $K \in \K_2$, we have
\[
c_{tr}(K) = \frac{ \area(K) + 2\area(P) }{\area(K)} = 1 + \frac{2\area(P) }{\area(K)}.
\]
By Theorem 3 of \cite{RS58}, the latter expression is maximal if, and only if $K$ is a convex quadrilateral, and by Theorem 1 of \cite{GHL13},
it is minimal if, and only if, $K$ is an ellipse. Thus the assertion readily follows.

Our next case is the left-hand side inequality of \ref{thm:main}.4.
Observe that
\begin{equation}\label{eq:mstarfirst}
c_{tr}^{m^*}(K) = \frac{4 \left( \area(K) + 2 \area(P) \right)}{\area(P')},
\end{equation}
where $P$ is a largest area inscribed, and $P'$ is a smallest area circumscribed parallelogram in $M = \frac{1}{2}(K-K)$.

As in the previous sections, if $c_{tr}^{m^*}(K)$ is minimal for some $K \in \K_2$, then, by \cite{C66}, we may assume
that $K$ is a Reuleaux triangle in its relative norm, and its area is $\area(K) = 2\area(M)-\frac{4}{3}\area(H)$,
where $H$ is a largest area affine-regular hexagon inscribed in $M$.
Thus, $\area(M) \geq \area(H)$ implies that
\begin{equation}
c_{tr}^{m^*}(K) \geq \frac{8 \left( \area(M) +  3\area(P) \right)}{3\area(P')}.
\end{equation}

On the other hand, we clearly have $\area(P) \geq \frac{1}{2} \area(P')$, where we have equality, for example, if $M$ is an affine-regular hexagon.
Furthermore, Corollary 5.1 of \cite{PT05} states that Gromov's mass* of any $o$-symmetric convex disk is at least three, with equality if, and only if $M$ is an affine-regular hexagon. This implies that $\area(M) \geq \frac{3}{4} \area(P')$, and thus, we obtain $c_{tr}^{m^*}(K) \geq 6$.
Here, we have equality if, and only if $M$ is an affine-regular hexagon, which immediately implies that $K$ is a triangle.

Finally, we prove the right-hand side of \ref{thm:main}.4. Similarly like in the previous sections, we may assume that $K = M$.
But then, clearly, $\area(M) \leq \area(P')$, $\area(P) \leq \area(P')$ and (\ref{eq:mstarfirst}) yields that $c_{tr}^{m^*}(K) \leq 12$.
Since in both inequalities equality is possible only if $M$ is a parallelogram, the assertion follows.

\section{Concluding remarks and open problems}\label{sec:remarks}

Our first question is to find the plane convex bodies $K$ for which the quantity $c_{tr}^{HT} (K)$ is maximal.

\begin{prob}
Prove or disprove such that if $c_{tr}^{HT} (K)$ is maximal for some $K \in \K_2$, then $K$ is an affine image of the body $M_0$ described in the Introduction.
\end{prob}

\begin{rem}\label{rem:centralsymmetrization}
For any $K \in K_n$ and direction $u \in \S^{n-1}$, let $d_u(K)$ denote the length of a maximal chord of $K$ in the direction $u$, and let
$K|u^\perp$ be the orthogonal projection of $K$ onto the hyperplane, through $o$, that is perpendicular to $u$.
Then the maximal volume of the convex hull of two intersecting translates of $K$ (that is, the numerator in the definition of $c_{tr}(K)$), is
\begin{equation}\label{eq:higherdimension}
\lambda_n(K) + \max \{ d_u(K) \lambda_{n-1}(K|u^\perp) : u \in \S^{n-1} \}.
\end{equation}
This observation can also be found in the proof of Theorem 1 of \cite{GHL13}.
Note that for any $u \in \S^{n-1}$, the central symmetral of $K|u^\perp$ is $\left(\frac{1}{2}(K-K)\right) |u^\perp$.
Thus, by the Brunn-Minkowski Inequality, the expression in (\ref{eq:higherdimension}) does not decrease under central symmetrization,
with equality if, and only if $K$ is centrally symmetric.
This yields that if $c_{tr}^{\tau}(K)$ is maximal for some $K \in \K_n$ for any $\tau \in \{ Bus, HT, m, m*\}$, then $K$ is centrally symmetric.
\end{rem}

\begin{rem}
By Remark~\ref{rem:centralsymmetrization}, to find the maximal value of $c_{tr}^{Bus}(K)$, it suffices to find the maximum of
$c_{tr}(K)$ over the family of $n$-dimensional centrally symmetric convex bodies. Thus, from Theorem 3 of \cite{RS58} it follows that
\[
c_{tr}^{Bus}(K) \leq n+1,
\]
with equality if, and only if $K$ is a centrally symmetric pseudo-double-pyramid in the sense of \cite{RS58}.
Similarly, by \cite{RS58} and \cite{GHL13}, over the family of $n$-dimensional $o$-symmetric convex bodies, we have
\[
c_{tr}^{Bus}(K) \geq 1+\frac{2v_{n-1}}{v_n},
\]
with equality if, and only if $K$ is an ellipsoid.
\end{rem}

\begin{prob}
For $n \geq 3$ and $\tau \in \{ HT, m, m*\}$, find the maximal values of $c_{tr}^\tau(K)$ over $\K_n$.
\end{prob}

\begin{prob}\label{prob:minimum}
For $n \geq 3$ and $\tau \in \{ Bus, HT, m, m*\}$, find the minimal values of $c_{tr}^\tau(K)$ over $\K_n$.
\end{prob}

When finding the minimal value of $c_{tr}^{Bus}(K)$ over $K \in \K_2$, we had to examine the smallest area convex disks of constant width two
in a fixed normed plane. Nevertheless, in $\Re^3$, even for the Euclidean norm, this question has been open for a long while (cf. \cite{KW11}).

Other problems arise if, instead of two translates of a convex body, we consider other families related to the body.
This was done also by Rogers and Shephard,
who, among other objects, studied the extrema of the volumes of \emph{differences bodies} or \emph{reflection bodies}.
We remark that a more general treatment of this type of questions can be found in \cite{GHL13} (cf. also \cite{gho} or \cite{gho2}).

Our problem applied to the case of difference bodies has already appeared in the literature in a different setting.
The Busemann volume of the difference body of $K$ is $2^n$ for any $K \in \K_n$.
For Holmes-Thompson volume, its value is a constant multiple of the volume product of the central symmetral of $K$,
and thus, its maximum is attained for ellipsoids, and the problem of finding its minimum leads to the famous Mahler Conjecture.
For Gromov's mass, we have
\[
\frac{4^n}{n!} \leq \vol_M^m(K-K) \leq 2^n v_n
\]
for every $K \in \K_n$ (cf. \cite{PT10}), and these inequalities are sharp.
For Gromov's mass*, we have
\[
\vol_M^{m*}(K-K) \leq 4^n
\]
and finding its minimum is also connected to the Mahler Conjecture (cf. \cite{PT10}).

Another possibility is to examine the reflection bodies of $K$, which are defined as
the convex hull of $K$ with one of its reflections about some point $x \in K$.

\begin{defi}
Let $K \in K_n$ and $M=\frac{1}{2}(K-K)$.
For $\tau \in \{ Bus, HT, m, m* \}$, set
\begin{equation}
c_p^\tau(K) = \max \{ \vol_M^\tau(\conv (K \cup 2x-K)): x \in K \}.
\end{equation}
\end{defi}

\begin{prob}
For $n \geq 2$ and $\tau \in \{ Bus, HT, m, m*\}$, find the minimal and the maximal values of $c_p^\tau(K)$ over $\K_n$.
\end{prob}

\noindent {\bf Acknowledgements.} The author is grateful to \'Akos G.Horv\'ath for the valuable conversations that they had on the topics covered in this paper,
and to Endre Makai, Jr. for letting him know about a theorem that enabled him to find the maximum for Holmes-Thompson volume in (\ref{thm:main}.2) of Theorem~\ref{thm:main}.

\end{document}